\documentclass[12pt]{article}
\usepackage[latin1]{inputenc}
\usepackage{amssymb}
\usepackage{amsmath}
\usepackage[english]{babel}
\usepackage{cases}
\usepackage{fancyhdr}
\usepackage{graphicx}
\usepackage{color}
\usepackage[colorlinks=true,linktocpage=true,linkcolor=blue,citecolor=red]{hyperref}

\usepackage[latin1]{inputenc}
\usepackage{color}
\usepackage{latexsym}
\usepackage{amssymb}
\usepackage{graphicx}
\usepackage{enumerate}
\usepackage{extarrows}
\newtheorem{Theorem}{Theorem}[part]
\newtheorem{Definition}{Definition}[part]
\newtheorem{Proposition}{Proposition}[part]
\newtheorem{Assumption}{Assumption}[part]
\newtheorem{Lemma}{Lemma}[part]
\newtheorem{Corollary}{Corollary}[part]
\newtheorem{Remark}{Remark}[part]

\newcommand{\vvertiii}[1]{{\vert\kern-0.25ex\vert\kern-0.25ex\vert #1
    \vert\kern-0.25ex\vert\kern-0.25ex\vert}}

\def \s{~~~}
\def \2{\vspace{2mm}}

\def \nn {\nonumber}

\def \Int{\displaystyle\int}

\def \R{\mathbb{R}}

\def \E{\mathbb{E}}
\def \F{\mathbb{F}}

\def \P{\mathbb{P}}
\def \H{\mathbb{H}}
\def \Q{\mathbb{Q}}

\def \S{\mathbb{S}}
\def \T{\mathbb{T}}

\def \Dc{{\cal D}}

\def \Fc{{\cal F}}

\def \Hc{{\cal H}}

\def \Pc{{\cal P}}

\def \Sc{{\cal S}}
\def \Tc{{\cal T}}

\def \eps{\varepsilon}

\def \ep{\hbox{ }\hfill$\Box$}

\def\t{\tau}
\def\th{\theta}
\def\ld{\lambda}

\def\be{\begin{eqnarray}}
\def\ee{\end{eqnarray}}
\newcommand{\ba}{\begin{array}}
\newcommand{\ea}{\end{array}}

\def\es{\mathop{\mbox{\rm ess\,sup }}}
\def \qq {\qquad}

\def\b*{\begin{eqnarray*}}
\def\e*{\end{eqnarray*}}
\def \nn{\nonumber }

\def\beqs{\begin{eqnarray*}}
\def\enqs{\end{eqnarray*}}
\def\beq{\begin{eqnarray}}
\def\enq{\end{eqnarray}}

\def \Yb {\overline{Y}}
\def \Zb {\overline{Z}}
\newcommand{\stopo}{\mathcal{T}_{0,T}}
\newcommand{\stops}{\mathcal{T}_{S,T}}
\DeclareMathOperator*{\esssup}{ess\,sup}

\begin{document}
\title{STRONG ENVELOPE AND STRONG SUPERMARTINGALE: \\APPLICATION TO REFLECTED BACKWARD STOCHASTIC DIFFERENTIAL EQUATION}
\author{Soufiane Aazizi \thanks{  Department of Mathematics, Faculty of Sciences Semlalia Cadi Ayyad University, B.P. 2390 Marrakesh, Morocco. \texttt{aazizi.soufiane@gmail.com}} \s \s Youssef Ouknine \thanks{ Department of Mathematics, Faculty of Sciences Semlalia Cadi Ayyad University, B.P. 2390 Marrakesh, Morocco. \sf ouknine@ucam.ac.ma}}


\maketitle

\begin{abstract} We provide several characterizations to identify Strong envelop (for bounded measurable process) and Strong super-martingale (for non-negative right upper semi-continuous process of the class $\Dc$).  As examples of application, we prove existence and uniqueness of reflected backward stochastic differential equation with lower barrier (RBSDB in short) in two cases: $i)$. the obstacle is a measurable bounded process; $ii)$. the obstacle is a right upper semicontinuous  optional process of class $\Dc$.

\end{abstract}

\vspace{13mm}

\noindent {\bf Key words~:} Strong envelope; Strong supermartingal; Snell envelope; Stochastic variational inequality (SVI).

\vspace{5mm}


\vspace{13mm}

\setcounter{equation}{1} \setcounter{Assumption}{1}
\setcounter{Theorem}{1} \setcounter{Proposition}{1}
\setcounter{Corollary}{1} \setcounter{Lemma}{1}
\setcounter{Definition}{1} \setcounter{Remark}{1}

\section{Introduction and preliminary proposition}
The Snell envelop is one of the fundamental concept used to solve the optimal stopping problem and reflected backward stochastic differential equation (RBSDE in short). But  recently, it has been pointed out by Peng and Xo \cite{PX05}, Essaky et al. \cite{EHO12} and Grigorova et al. \cite{GIOOQ15} that a more generalized versions of Snell envelop solve RBSDE with irregular barrier. Hence the interest to study these generalizations and provide theme applications.\\
To the best of our knowledge, the most generalizations of Snell envelop of c�dl�g process are: Strong supermartingale of nonnegatif l�dl�g optional process of class $\Dc$ introduced by Mertens \cite{M72} and Strong envelop of measurable bounded process introduced by Stettner  and Zabczyka \cite{SZ81}.\\

The present work is divided in three parts. In the first part, we provide an important characterization of Strong supermartingale, which is defined according to Mertens \cite{M72} as the smallest supermartingale majoring a given l�dl�g non-negative optional process. Here, in our setting, we are concerned by the framework of Kobylanski and Quenez \cite{KQ12} of non-negative right upper semicontinous (USC in short) process $Z$, and we will show that the process $Y$ is the Strong supermartingale of $Z$, if and only if the following inequality and Skorohod conditions hold true:
    \begin{equation*}
\left\{
\begin{array}{lll}
    &&  Z \geq Y,   \text{ up to an evanescent set }\\
         &&\Int_0^{T} 1_{\{Z_t>Y_t\}}dA^c_s=0, \mbox{ a.s., }\\
         &&\bigtriangleup C_t=\bigtriangleup C_t1_{\{Z_t=Y_t\}}=(Z_t-Z^+_t)1_{\{Z_t=Y_t\}}\mbox{ a.s., }\\
         &&\bigtriangleup A^d_t=\bigtriangleup A^d_t1_{\{Z_{t^-}=\bar{Y}_t\}}\mbox{ a.s., }
         \end{array}
         \right.\end{equation*}
where $A$ and $C$ are  the increasing processes appearing in Mertens decomposition of $Y$, and $\bar{Y}_t=\limsup_{s\rightarrow t, s< t} Y_s, \qq 0 \le t \le T, \mbox{ a.s.}$

In the second part, we give two characterization of Strong envelop defined as the smallest supermartingale  majoring a given measurable bounded process, using the general Skorohod condition. Indeed, let $X$ be a measurable bounded process, $X^*$ be a measurable c�dl�g process such that $\E\sup_{0 \le t \le T }| X^*_t|^2<\infty$, and $U$ a Strong envelop, such that $X\le X^* \le U $ a.s., a.e.. Then we show that $U$ is the strong envelop of $X$ if and only if
    \begin{eqnarray*}
           \int_0^{T} (U_{s^-}-  X^*_{s^-})dA_s=0, \mbox{ a.s.. }
           \end{eqnarray*}
The second characterization of Strong envelop claims that the unique solution of stochastic variational inequality (SVI in short) is given by Strong envelop $U$ of the process $X$. In fact, let define a convex subset of the space of c�dl�g process $V$ such that $\E\sup_{0 \le t \le T }| V_t|^2<\infty$, and taking the following form
\begin{eqnarray*}
\mathcal{K}=\{V \in \S^2: \s V \geq X \s dt\otimes d\P-\mbox{a.e. and }  V_T=0  \mbox{ a.s. }\},
\end{eqnarray*}
and let for any $t\in [0, T]$, $X_t \le K$ a.s., where $K$ is a positive constant. Then, we will prove that the unique solution of the following SVI
  \begin{eqnarray*}\mathbb{E} \left[\int_{\tau_1}^{\tau_2} (U_{s-}  - V_{s-})dU_s/\mathcal{F}_{\tau_1} \right] \geq 0 \s \mbox{ a.s. }
\end{eqnarray*}
is $U$ the Strong envelop of $X$. Hence, we provide new estimate for the increment of the predictable component of the Strong envelope on an arbitrary stochastic interval through the supremum of the increments of the given process X on the same time interval

The third and last part of this paper, deals with two applications of the above characterization on RBSDE. We will extend the recent work of Grigorova et al. \cite{GIOOQ15} and we will establish existence and uniqueness of a RBSDE with jump taking the form
\begin{equation}
\left\{
\begin{array}{lll}
 Y_\tau=\xi_T+\Int_{\tau}^T f(t,Y_t,  Z_t)dt-\Int_{\tau}^T  Z_t dW_t -\int_{\tau}^T \int_U k_s(u)\tilde{N}(ds,du) +A_T-A_\tau+C_{T-} -C_{\tau-} \nn\\
 \text{ a.s. for all }\tau\in\stopo, \\
  Y_t \geq \xi_t   \text{ for all } t\in[0,T]  \text{ a.s.,} \\
 A \text{ is a nondecreasing right-continuous predictable process
with } A_0= 0, E(A_T)<\infty \text{ and } \nonumber\\
\text{ such that } \nonumber\\
 \Int_0^T {\bf 1}_{\{Y_t > \xi_t\}} dA^c_t = 0 \text{ a.s. and } \; (Y_{\tau-}-\bar{\xi}_\tau)(A^d_{\tau}-A^d_{\tau-})=0 \text{ a.s. for all (predictable) }\tau\in\stopo, \\
 \text{where } \bar{\xi}_t:=  \limsup_{s \uparrow t, s< t} \xi_s. \nn\\
 C \text{ is a nondecreasing right-continuous adapted purely discontinuous process }\\ \\
 \text{ with } C_{0-}= 0, E(C_T)<\infty  \text{ and such that }
(Y_{\tau}-\xi_{\tau})(C_{\tau}-C_{\tau-})=0 \text{ a.s. for all }\tau\in\stopo.
\end{array}
\right.
\end{equation}
where the obstacle $\xi$ is an optional right USC process of class $\Dc$. We also make the link between our BSDE and optimal stopping with $f$-conditional expectation. We then characterize the value of the problem in terms of the unique solution of RBSDE associated to $\xi$ and the driver $f$.

In the second application, we are concerned by proving existence and uniqueness of a RBSDE where the obstacle is not regular, in the spirit of the problem treated in the paper of Peng and Xu \cite{PX05} and also in Essakey et al. \cite{EHO12}.
In our setting, we will show in Theorem \ref{RBSDEone} that $Y$ is the minimal solution of the following RBSDE with lower barrier $L$
\begin{equation}
\left\{
\begin{array}{lll}
Y_t = -\Int_0^tf(s, Y_s, Z_s)ds+\Int_t^T dK_s- \int_t^T Z_sdB_s; \\
Y\geq L \s \mbox{ a.s. a.e.};\\
\Int_0^{T}(Y_{s^-} - L^*_{s^-})dK_s=0, \s \forall L^* \in D_\mathcal{F}\mbox{ such that } L \leq L^* \leq Y \s \mbox{a.s. a.e.}
\end{array}
\right.
\end{equation}
Furthermore, we show that $\Big(Y_t + \Int_0^tf(s)ds\Big)_{0 \le t \le T }$ is the Strong envelope of $\Big(L_t + \Int_0^tf(s)ds\Big)_{0 \le t \le T }$.\\


In section \ref{CharactSSM} we characterize the solution of Strong super-martingale using Mertens decomposition and Skorohod condition. In section \ref{Charact}, we characterise the solution of Strong envelope using Doob-Meyer decomposition and the generalized Skorohod condition. In Theorem  \ref{ThSVI}, we prove that the unique solution of the SVI is given by the Strong envelope $U$ of the given measurable bounded process $X$.  In section \ref{Section Estimate} we formulate and prove the estimate for the increment of the predictable process component of the Strong envelope on an arbitrary stochastic interval through the supremum of the increment of the given process $X$ on the same interval. We also provide an estimate on the square bracket of the difference of two Strong envelope. Finally, in section \ref{Appli}, we show how we can use the notion of Strong envelope and Strong super-martingale to prove existence and uniqueness of reflected backward stochastic differential equation with a lower barrier.

\section{Strong super-martingale}\label{CharactSSM}
The notion of Strong super-martingale was developed by Mertens \cite{M72} to define the smallest supermartingale of an optional process $Y$, this result was  generalized by Elkaroui \cite{E81} to  prove the existence of an optimal stopping time when the reward is given by an upper semi-continuous non negative process of class $\Dc$. Recently, Kobylanski and Quenez \cite{KQ12} resolved the problem under week assumption in term of integrability and regularity of the reward family, supposing that that reward family is upper semicontinous along stopping times in expectation.\\
We recall some important properties proved by Maingueneau \cite{M78}, which will be used later in the paper.
\begin{Proposition}\label{Proposition Snell envelop}
Let  $Y=(Y_t)_{ 0 \le t \le T}$ be an optional positive process, that belongs to the class $\Dc$. Then there exist a unique, optional process $Z=(Z_t)_{0 \le t \le T}$  such that $Z$ is the smallest strong super-martingale which dominate $Y$, the process $Z$ is called the Snell envelope of $Y$ and it has the following properties:
\begin{enumerate}
    \item For any stopping time $\th$, we have
    \be
    Z_\th= SN(Y)_\theta= \es_{\t\in\T_\th}\E[Y_\t/\Fc_\th],\qq Z_{T}=Y_{T}.
    \ee
\item Mertens decomposition: There exist a uniform c�dl�g martingale $M$  and unique predictable right continuous non decreasing process $(A_t)$ with $A_0=0$ and $\E[A_T]<\infty$ and a unique right continuous adapted  non decreasing process $(C_t)$, which is purely discontinuous with $C_0=0$ and $\E[C_T]<\infty$, such that
    \be
Z_t=M_t-A_t-C_{t-}, \qq 0\le t \le T, \mbox{ a.s.}
\ee
We have $ \bigtriangleup C_t= Z_t - Z_{t+}$. And recall that the process $(A_t)$ admits the following unique decomposition: $A_t=A^c_t+A^d_t$, where $(A^c_t)$ is the continuous part of $(A_t)$ and $(A^d_t)$ is its purely discontinuous part.

\item A stopping time $\th$ is said to be optimal if we have $Z_\th=Y_\th$ a.s. and $Z_{t\wedge \th}$ is a martingale;
\item Let $\ld \in [0, 1[$, and denote $J^\ld=\{(\omega, t) / Y_t(\omega) > \ld Z_t(\omega)\}$  and \\
\b*
D^\ld_t(\omega)=inf\{s\ge t / (\omega, s) \in J^\ld\}\wedge T, \qq (\mbox{we write } D_0^\ld=D^\ld);
\e*
Then, for any stopping time $\th$, we have a.s.
\be
Z_\th=\E[Z_{D^\ld_\th} / \Fc_\th];
\ee
\item    We have $\ld Z_{D^\ld_\th} \le Y_{D^\ld_\th}\vee Y^+_{D^\ld_\th}$;
\end{enumerate}
\end{Proposition}
Note that according to the result of General Theory of processes in El Karoui \cite{E80}, for each adapted process $(Y_t)$, there exists a predictable process $(\bar{Y}_t)$ such that
\b*
\bar{Y}_t=\limsup_{s\rightarrow t, s< t} Y_s, \qq 0 \le t \le T, \mbox{ a.s.}
\e*
By the following proposition, we characterise Strong super-martingale of a given process using Skorohod condition
\begin{Proposition} \label{PropSE(X)}
Let $Z$ be a Strong super-martingale and $Y$ a right upper semi-continuous process (USC for short) process of the class $\Dc$, we then have the following equivalence:
\begin{enumerate}
  \item $Z$ is the Strong Snell envelope  of $Y$.
  \item \begin{description}
          \item[a.)] $Z \geq Y$,   up to an evanescent set;
          \item[b.)] The following Skorohod condition  holds true:
           \begin{eqnarray}
           \int_0^{T} 1_{\{Z_t>Y_t\}}dA^c_s=0, \mbox{ a.s., }
           \end{eqnarray}
           \item[c.)] $\bigtriangleup C_t=\bigtriangleup C_t1_{\{Z_t=Y_t\}}=(Z_t-Z^+_t)1_{\{Z_t=Y_t\}}$ a.s.
           \item[d.)]  $\bigtriangleup A^d_t=\bigtriangleup A^d_t1_{\{Z_{t^-}=\bar{Y}_t\}}$ a.s.
                     \end{description}
        where $A$ is the increasing process appearing in Doob-Meyer decomposition of $U$.
\end{enumerate}
\end{Proposition}
{\bf Proof.}\\
To prove $ 1 \Longrightarrow 2$. we send the reader to Proposition B.11 in \cite{KQ12}\2 \\
To prove $ 2 \Longrightarrow 1$. Let us consider the stopping time for each $\omega$ by $\t^\ld_\th(\omega)=\inf\{t\ge\th(\omega)/ Y_\t(\omega) >\lambda Z_\t(\omega)\}$. Hence, from Proposition B.5 in \cite{KQ12} we have
\be\lambda Z_{\t^\ld_\th} \le Y_{\t^\ld_\th} \mbox{ a.s.}\ee
Taking conditional expectation in both hand side in the above inequality combined with the fourth point of Proposition \ref{Proposition Snell envelop},  we get
\be\lambda Z_\theta &\le& \E[Y_{\t^\ld_\th}/\Fc_\th] \\
                &\le& SN(Y)_\th\ee
By sending $\ld$ to $1$, we obtain $Z_\th\le SN(Y)_\th$, then $ Z\le SN(Y)$ up to an evanescent set (u.e.s. in short). From other side, we have $Y\le Z$ u.e.s. which implies that $SN(Y)\le Z$ u.e.s. since $SN(Y)$ is the smallest super-martingale majoring $Y$. Hence, we deduce from the last inequalities that $Z=SN(Y)$ u.e.s.  which conclude the proof.\ep

\section{Strong envelope}\label{Charact}
Let $T$ be a fixed positive real number, Let $(\Omega, \Fc, (\Fc_t)_{0 \le t \le T }, \P)$ be a probability space, satisfying the usual conditions, where $(\Fc_t)_{0 \le t \le T }$ is an increasing right-continuous family of complete $\sigma$-fields. We introduce the following notation:\\
$[\mathfrak{B}]$\quad: the class of bounded progressive measurable process $X$.\\
$\S^2$ \,\quad: the Banach space of all c�dl�g process $Y$ such that
\beq
\| Y \|_{\S^2}:=\Big(\E\Big[\sup_{0 \le t \le T }| Y_t|^2\Big]\Big)^{\frac{1}{2}} < \infty,
\enq
$H^p\quad:$ $1 \le p < \infty$, the space of semimartingales $Y$ (see \cite{DM82}) such that
\beq
\| Y \|_{H^p}=\inf_{X=M+A}\left(\left\|[ M]_T^\frac12+ \int_0^{T} |dA_s|\right\|_{L^p}\right)< \infty,
\enq
where the infimum is taken over all decompositions of $Y$ into local martingale and the process of bounded variation. \\
Let $(X_t)_{0 \le t \le T }$ be a bounded progressively measurable process.
A right-continuous super-martingale $U:=(U_t)_{0 \le t \le T }$ is called the Strong envelope of $X:=(X_t)_{0 \le t \le T }$ if it is the smallest right-continuous, non-negative super-martingale  such that $ U \geq X,$ $dt\otimes d\P $-a.e.; i.e., if $(\bar{U}_t)_{0 \le t \le T } $ is another super-martingale such that $\bar{U}\geq X,$ $dt\otimes d\P $-a.e., then $\bar{U}_t\geq U_t,$ a.s. for any $t \in [0, T]  $ and we write $U=SE(X)$.\\
The problem of Strong envelope is formulated as the following: For an arbitrary $\beta >0 $  find a right-continuous process $U^\beta$ such that for all $t \in [0, T]$
        \beq
        \label{Ubeta}U_t^\beta:= \beta \mathbb{E} \left(\int_t^{T}(X_s - U^\beta_s)^+ds / \mathcal{F}_t\right) \s\P-\mbox{a.e.}
        \enq
If we assume tat $X$  is a bounded, progressively measurable process, then according to Theorem 9 in \cite{Z83}, there is a solution of (\ref{Ubeta}) and is unique, up to indistinguishable process. It increases with $ \beta \nearrow +\infty$ and the limit process \begin{eqnarray}
U_t:= \lim_{\beta \nearrow +\infty } U_t^\beta, \s \s 0 \le t \le T ,
\end{eqnarray}
is the Strong envelope of $X$.
\subsection{Characterization of Strong envelope}
In what follow, we give the main tools to characterize the Strong envelope of a given process. Indeed, we provide in the following proposition the first characterization of Strong envelope. 
\begin{Proposition} \label{PropSE(X)}
Let $U$ be a super-martingale and $X$ a process of the class $[\mathfrak{B}]$, then we have the following equivalence:
\begin{enumerate}
  \item $U$ is the Strong envelope  of $X$.
  \item \begin{description}
          \item[a.)] $U \geq X$,  a.s. a.e., and $U_T=0$ a.s.;
          \item[b.)] The following (generalized) Skorohod condition (cf. \cite{S65}) holds true:
           \begin{eqnarray}
           \int_0^{T} (U_{s^-}-  X^*_{s^-})dA_s=0, \label{GSC}  \mbox{ a.s., } \forall X^* \in \S^2\mbox{ s.t. }  U \geq X^* \geq X \s\mbox{a.s., a.e.}
           \end{eqnarray}
        \end{description}
        where $A$ is the increasing process appearing in Doob-Meyer decomposition of $U$.
\end{enumerate}
Moreover, if $X^*_T=0$ a.s. then, we have (a.s.)
\begin{eqnarray}
U_t=\mathbb{E}\left[U_{\tau^{\eps*}_t} /\mathcal{F}_t\right].
\end{eqnarray}
for any $\eps >0$ and $t\le T$, where the stopping time $\tau^{\eps*}_t$ is given by
\beq \label{tau}
\tau^{\eps*}_t=\inf\{s \geq t,  \mbox{ s.t.} \s  X^*_s \geq U_s - \eps \}\wedge T,
\enq
\end{Proposition}
{\bf Proof.}\\
We prove $ 1 \Longrightarrow 2$. \2 \\
Let $U=SE(X)$, then (2.a) is obvious. \2 \\
Let $X^*$ be a c�dl�g process such that $U \geq X^*\geq X$ a.s. a.e., then $U_- \geq X^*_-$ a.s. and $U \geq SN(X^*) \geq X$ a.s., a.e., where $SN(X^*)$ is the Snell envelope of $X^*$. However, $U$ is the smallest super-martingale that majors $X$ then $U \leq SN(X^*)$ a.s., combining this with inverse inequality we have $U = SN(X^*)$. By the uniqueness of Doob-Meyer decomposition of $U_t= M_t - A_t$ and $SN(X^*_t)= M^*_t - A^*_t$, we have $M_t^*=M_t$ and $A^*_t=A_t$, which leads to
\begin{eqnarray*}
\int_0^. 1_{\{{U_t}_- > {X^*_t}_- \}}dA_s = \int_0^. 1_{\{{SN(X^*)_t}_- > {X^*_t}_- \}}dA_s^*=0.
\end{eqnarray*}
To prove the inverse inequality, we consider a c�dl�g process $X^*$  such that $X \leq X^* \leq U \wedge SE(X)$, and let the stopping time $\tau^{\eps*}_t$ be such that
\begin{eqnarray*}
\tau^{\eps*}_t=\inf\{s \geq t,  \mbox{ s.t.} \s  X^*_s \geq U_s - \eps \}.
\end{eqnarray*}
Then
\begin{eqnarray*}
\mathbb{E}\left[U_{\tau^{\eps*}_t} - U_t /\mathcal{F}_t\right] &=& - \mathbb{E}\left[\int 1_{t < s \leq \tau^{\eps*}_t} dA_s /\mathcal{F}_t\right]\leq 0.
\end{eqnarray*}
Notice that $$(t, \tau^{\eps*}_t] =\{(s, \omega) / \s t< s \le \tau^{\eps*}_t(\omega)\}~\subset ~\{(s, \omega)/ \s s> t \mbox{ and } X^*_{s-}< U_{s-}  \},$$ which leads to
\begin{eqnarray*}
 \mathbb{E}\left[\int 1_{t \leq s < \tau^{\eps*}_t} dA_s /\mathcal{F}_t\right]\le  \mathbb{E}\left[1_{U_{s-} > X^*_{s-}}dA_s /\mathcal{F}_t\right]= 0.
\end{eqnarray*}
Combining this with the latest above inequality, we conclude the martingale property of $U$ and we have
\begin{eqnarray}\label{tauStart}
U_t=\mathbb{E}\left[U_{\tau^{\eps*}_t} /\mathcal{F}_t\right].
\end{eqnarray}
From other side, $U_{\tau^{\eps*}_t} \leq X^*_{\tau^\eps_t} + \eps $
taking the conditional expectation in both side
\begin{eqnarray*}
U_t &\leq& \mathbb{E}\left[X^*_{\tau^{\eps*}_t}  /\mathcal{F}_t\right] + \eps \\
    &\leq& SN(X^*_t) + \eps.
\end{eqnarray*}
By sending $\eps$ to $0$, we get $SN(X^*)= SE(X)=U$. \ep
\begin{Remark}
If we assume furthermore that $X$ is a c�dl�g process, then $U$ is the snell envelope of $X$ and the Skorohod condition $(\ref{GSC})$ becomes
\beqs
           \int_0^T (U_{s^-}-  X_{s^-})dA_s=0.
           \enqs
\end{Remark}
\begin{Remark}\label{U in S2}
If there exists a process $X^* \in \S^2$ such that $X\le X^* \le U:=SE(X)$ a.s. a.e. Then $U\in \S^2$.
\end{Remark}
The second characterization of Strong envelope states that the unique solution of the stochastic variational inequality (SVI in short) is given by the Strong envelope $U$  of the given process $X$. Indeed, let formulate the SVI related to optimal stopping problem, and Let the process $X$ be in $[\mathfrak{B}]$ and $ \mathcal{K}$ be a convex subset of the space $\S^2$, taking the following form
\begin{eqnarray}
\mathcal{K}=\{V \in \S^2: \s V \geq X \s dt\otimes d\P-\mbox{a.e. and }  V_T=0  \mbox{ a.s. }\}.
\end{eqnarray}
The problem of stochastic variational inequality associated to optimal stopping time, consists to find an element $U \in \mathcal{K} \cap H^2$ such that for any element $V \in \mathcal{K}$, any pair of stopping times $(\tau_1, \tau_2)$ where $0\leq \tau_1 \leq  \tau_2 $,  the following inequality holds
\begin{eqnarray}\label{SVI}
\mathbb{E} \left[\int_{\tau_1}^{\tau_2} (U_{s-}  - V_{s-})dU_s/\mathcal{F}_{\tau_1} \right] \geq 0 \s \mbox{ a.s. }
\end{eqnarray}
The conditional expectation in the above inequality is well defined, in fact, using Emery's inequality (See \cite{P90}, Chapter V, Theorem 3) we obtain
\begin{eqnarray*}
\left\|\int_{\tau_1}^{\tau_2} (U_{s-}  - V_{s-})dU_s\right\|_{H^1} &=&\left\|\int_0^T1_{\tau_1<s\le \tau_2} (U_{s-}  - V_{s-})dU_s\right\|_{H^1}  \\&\le &\left\|U-V\right\|_{\S^2}  \cdot \left\|U\right\|_{H^2} < \infty.
\end{eqnarray*}
We are now in the position to provide  the second characterization of the Strong envelope:
\begin{Theorem}\label{ThSVI} Let $X$ be a bounded progressive measurable process such that for any $t\in [0, T]$ we have $X_t\le K$ a.s. where $K$ is a positive constant. Then, there exists a unique solution $U$ of the SVI \eqref{SVI}, which is the Strong envelope of $X$.
\end{Theorem}
{\bf Proof.} Let us check that $U=SE(X)$ is a solution of above SVI. \\
we have $ U \in \mathcal{K} \cap H^2$. In fact, since $X$ is bounded by $K$ it follows that $U$ is also bounded by $K$ ($c.f.$ definition \eqref{Ubeta}).
From other side, by Proposition \ref{PropSE(zn)}, we know that $ U_t = M_t -A_t$, which leads to
\begin{eqnarray*}
\E(A_T-A_t/ \mathcal{F}_t)=U_t&\le& K.
\end{eqnarray*}
Hence, \begin{eqnarray}\label{AinLp}
\left\|A_T|\right\|_{L^p}&\le &K.
\end{eqnarray}
Combining this with the following consequence
\begin{eqnarray*}
\sup_{0 \le t \le T }|M_t| &\le& \sup_{0 \le t \le T }|U_t|+A_T.
\end{eqnarray*}
Finally,
\begin{eqnarray}\label{MinLp}
\left\|\sup_{0 \le t \le T }|M_t|\right\|_{L^p} &\le& 2K.
\end{eqnarray}
From \eqref{AinLp} and \eqref{MinLp} we deduce that $\|U\|_{H^2} < \infty$ and since $U\ge X$-a.s. a.e. and $U_T=0$-a.s, then we have $U\in \mathcal{K}\cap H^2$.\\
Taking  $V \in \mathcal{K}$, and consider the stochastic integral
$$ \int_{\tau_1}^{\tau_2}[U_{s-} - V_{s-} ]dU_s.$$
Since this  stochastic integral belongs to the class $H^1$, then its martingale part vanishes after taking conditional expectation. From other side, since we have $U \leq SN(V)$, it follows by adding and subtracting $SN(V)$ in the above stochastic integral that
\beqs
\mathbb{E} \left[\int_{\tau_1}^{\tau_2}[U_{s-} - V_{s-} ]dU_s/\mathcal{F}_{\tau_1} \right]   &=&- \mathbb{E} \left[\int_{\tau_1}^{\tau_2}[SN(V)_{s-} -U_{s-} ]dA_s/\mathcal{F}_{\tau_1} \right]\\
&&+ \mathbb{E} \left[\int_{\tau_1}^{\tau_2}[SN(V)_{s-} - V_{s-} ]dA_s/\mathcal{F}_{\tau_1} \right] \\
&\geq& 0.
\enqs
To conclude the proof, we need the uniqueness of the solution. Let $U$, $U'$ $\in \mathcal{K} \cap H^2$ two solutions of SVI (\ref{SVI}). Then
\begin{eqnarray*}
\mathbb{E} \left[\int_{\tau_1}^{\tau_2} (U_{s-}  - U'_{s-})d(U_s - U'_s)/\mathcal{F}_{\tau_1} \right] \geq 0.
\end{eqnarray*}
Taking $\tau_1=t$ , $\tau_2=T$, and using It� formula on $(U-U')^2$ leads to
\begin{eqnarray*}
 - (U_{t}-U'_{t})^2 &=& 2 \int_{t}^{T}(U_{s-}-U'_{s-})d(U_s-U'_s) + [U-U']_{T} - [U-U']_{t}.
\end{eqnarray*}
Taking the conditional expectation, we obtain
\begin{eqnarray*}
 - (U_{t}-U'_{t})^2 &=& 2 \mathbb{E}\left[\int_{t}^{T}(U_{s-}-U'_{s-})d(U_s-U'_s)/ \mathcal{F}_{\tau_1} \right] \\
 &&+ \mathbb{E}\left[[U-U']_{T} - [U-U']_{t}/ \mathcal{F}_{\tau_1} \right] \ge 0.
\end{eqnarray*}
Hence, for all $t\in [0, T]$, $U_t=U'_t $ a.s.
\ep \\
In the following proposition, we give the main convergence property of Strong envelope.
\begin{Proposition}\label{PropSE(zn)}
The process $U:=SE(X)$  enjoys the following properties:
\begin{enumerate}[(i)]
\item   The Doob-Meyer decomposition of the  super-martingale $U$ implies the existence of a martingale $(M_t)_{0 \le t \le T }$ and a non-decreasing processes $(A_t)_{0 \le t \le T }$ which is and predictable such that,
    $
    U_t=M_t-A_t, \, 0\le t\le T.$
    \item If $(X^n)_{n\geq 0}$ and $X$ are a bounded progressively measurable processes such that the sequence $(X^n)_{n\geq 0}$ is increasing and converges pointwisely to $X$, then $SE(X^n)_{n\geq 0}$ is also increasing and  converges pointwisely to  $SE(X)$.
\end{enumerate}
\end{Proposition}
{\bf Proof. }
The first assertion is obvious, we then prove $(ii)$.\\
Since $X^n$ converges increasingly to $X$, it follows that for all $0 \le t \le T $, $SE(X^n) \leq SE(X)$-a.s.. \\
Furthermore, \beq \label{v<SE(X)} \lim_{n \rightarrow \infty}SE(X^n) &\leq& SE(X) \mbox{ a.s.}\enq
Then $\lim_{n \rightarrow \infty}SE(X^n)$ is c�dl�g super-martingale (see e.g., Delacherie and Meyer \cite{DM82} p. 86). From other side, $X^n \leq SE(X^n)$ a.s., a.e., implies that $dt\otimes d\mathbb{P}$-a.e.  $X \leq \lim_{n \rightarrow \infty}SE(X^n)$. However, since $SE(X)$ is the smallest super-martingale which dominate $X$ then  $SE(X) \le \lim_{n \rightarrow \infty  }SE(X^n) $-a.s., together with (\ref{v<SE(X)}) we get the required result.\ep

\subsection{A priori estimate} \label{Section Estimate}
In this section, we provide an interesting a priori estimate, mainly based on the above characterization of Srong envelope.
On an arbitrary stochastic interval $(\sigma_1, \sigma_2]$, where $\sigma_1$ and $\sigma_2$ are two stopping times, we are going to formulate and prove a priori estimate for the increment of the predictable process $A$ component of Strong envelope. We adapt the proof of \cite{DDS03} to our setting
\begin{Theorem} \label{k-k}
Let $X$ be a bounded progressively measurable process, and $U$ its Strong envelope. Then, for any c�dl�g process $X^*$ such that
\beqs X \leq X^* \leq  U, \s \mbox{ a.s., a.e.}, \enqs
we have the following inequality for the increasing process $A$ component of the Strong envelope $U$:
\beq\label{A1-A2 le eps}
\mathbb{E}\Big( A_{\sigma_2} - A_{\sigma_1} / F_{\sigma_1} \Big) \leq \mathbb{E} \Big( X^*_{\tau^{\eps *}_{\sigma_1}\wedge \sigma_2 } - X^*_{ \sigma_2 }/ F_{\sigma_1} \Big) + \eps,
\enq
for arbitrary $\eps>0$, where $ \tau^{\eps *}_{\sigma_1}$ is defined in  \eqref{tau}. Moreover,  for any $p \geq 1$:
\beq\label{A1-A2 le p}\|A_{\sigma_2} - A_{\sigma_1}   \|_{L^p} \leq p \Big\|\sup_{\sigma_1 \leq s \leq \sigma_2} |X^*_{s} - X^*_{ \sigma_2 }| \Big\|_{L^p}.   \enq
\end{Theorem}
{\bf Proof.}
Denote $$\sigma(u)=(\sigma_1+u)\wedge \sigma_2, \s u\ge 0, $$
with $\sigma(0)=\sigma_1 $ and $\sigma(T)=\sigma_2 $. It follows that for $0 \le u \le t \le T$:
$$\E(A_{\sigma(t)}-A_{\sigma(u)}/\mathcal{F}_{\sigma(u)})= \E(U_{\sigma(u)}-U_{\sigma(t)}/\mathcal{F}_{\sigma(u)}).$$
Combining this with \eqref{tauStart}, we obtain
$$U_{\sigma(u)}=\E(U_{\tau^{\eps *}_{\sigma(u)}}/\mathcal{F}_{\sigma(u)}).$$
And by super-martingale property od $U$ we have
$$U_{\sigma(u)}=\E(U_{\tau^{\eps *}_{\sigma(u)}\wedge \sigma(t)}/\mathcal{F}_{\sigma(u)}).$$
Hence
\begin{eqnarray*}
\E(A_{\sigma(t)}-A_{\sigma(u)}/\mathcal{F}_{\sigma(u)})&=&\E(U_{\tau^{\eps *}_{\sigma(u)}\wedge\sigma(t)}-U_{\sigma(t)}/\mathcal{F}_{\sigma(u)})\\
&=&\E(1_{\tau^{\eps *}_{\sigma(u)}<\sigma(t)}\cdot (U_{\tau^{\eps *}_{\sigma(u)}}-U_{\sigma(t)})/\mathcal{F}_{\sigma(u)})\\
&\le&\E(1_{\tau^{\eps *}_{\sigma(u)}<\sigma(t)}\cdot (X^*_{\tau^{\eps *}_{\sigma(u)}}+\eps-U_{\sigma(t)})/\mathcal{F}_{\sigma(u)})\\
&\le&\E(X^*_{\tau^{\eps *}_{\sigma(u)}\wedge\sigma(t)}-X^*_{\sigma(t)}/\mathcal{F}_{\sigma(u)})+\eps.
\end{eqnarray*}
Taking $u=0$, $t=T$, we conclude estimate \eqref{A1-A2 le eps}. Taking $t=T$ in the above inequality and having in mind the following
$$ |X^*_{\tau^{\eps *}_{\sigma(u)}\wedge\sigma_2}-X^*_{\sigma_2}| \le \sup_{\sigma_1 \le s \le \sigma_2}|X^*_s-X^*_{\sigma_2}|, $$
we conclude by sending $\eps$ to $0$ that
$$\E(A_{\sigma_2}-A_{\sigma_1}-(A_{\sigma(u)}-A_{\sigma_1})/\mathcal{F}_{\sigma(u)})=\E\left(\sup_{\sigma_1 \le s \le \sigma_2}|X^*_s-X^*_{\sigma_2}|/\mathcal{F}_{\sigma(u)}\right).$$
Applying Garsia's inequality to the non-decreasing process $\hat{A}_u=A_{\sigma(u)}-A_{\sigma_1} $ we deduce the estimate \eqref{A1-A2 le p}.\ep
\begin{Corollary}
Let $U^i$, i=1, 2, be the Strong envelope of $X^i$ with $X^i \in [\mathfrak{B}]$.  Then for any arbitrary stochastic interval $(\tau_1, \tau_2]$ the following estimate hold:
\beqs &&\E\Big((U^2_{\tau_1} - U^1_{\tau_1})^2 +[U^2 - U^1]_{\tau_2} - [U^2 - U^1]_{\tau_1} \Big) \\
&& \leq 4 \Big\| \sup_{\tau_1 \leq t \leq \tau_2} |X^2_t -X^1_t| \Big\|_{L^2} \cdot \left(\Big\| \sup_{\tau_1 \leq t \leq \tau_2} |X^1_t -X^1_{\tau_2}| \Big\|_{L^2} + \Big\| \sup_{\tau_1 \leq t \leq \tau_2} |X^2_t -X^2_{\tau_2}| \Big\|_{L^2}    \right)\\
&&\quad+ \E(U^2_{\tau_2} - U^1_{\tau_2})^2 .
\enqs
In particular:
\beqs \E[U^2 - U^1]_T  &\leq& 4 \Big\| \sup_{0 \le t \le T } |X^2_t -X^1_t| \Big\|_{L^2} \\
&&\times \left(\Big\| \sup_{0 \le t \le T }|X^1_t -X^1_T| \Big\|_{L^2} + \Big\| \sup_{0 \le t \le T } |X^2_t -X^2_T| \Big\|_{L^2}    \right)+ \E(X^2_T -X^1_T)^2.
\enqs
\end{Corollary}
{\bf Proof.} The proof is similar to Theorem 2.3 in \cite{DDS03} using result of stochastic variational inequality. \ep

\section{Applications on Backward stochastic differential equations}\label{Appli}
\subsection{Reflected BSDE with right USC barrier}
This application extend the framework of \cite{GIOOQ15} to BSDEs  with jumps, discussed in \cite[Section 6]{GIOOQ15}.\\
Let $(\Omega, \Fc, \P)$  be a probability space equipped with a one-dimensional Brownian motion $B$ and an independent Poisson random measure $N(dt, du)$ defined on $(U, \mathcal{U})$ whose compensator is the measure $dt\otimes\mu(du)$, where $(U, \mathcal{U})$ is a measurable space equipped with a $\sigma$-finite positive measure $\mu$. Let $\tilde{N}$ be its compensated measure. The filtration  $\F=\{\Fc_t, t\ge 0\}$ corresponds to the complete natural filtration associated with  $B$ and $N$. Fot $t\in [0, T]$, we denote $\Tc_{0,T}$ the set of stopping times $\t$ such that $\P(t\le \t \le T)=1$. More generally, for a given stopping time $\nu \in \Tc_{0, T}$, we denote $\Tc_{\nu, T}$ the set of stopping times $\t$ such that $\P(\nu\le \t \le T)=1$.\\
We mean by $\mathcal{P}$ the $\sigma$-algebra of predictable set $\Omega\times [0, T]$.\\
We use the following notation
\begin{itemize}
\item
 $L^2({\Fc}_T)$  is the set of random variables which are  $\Fc_T$-measurable and square integrable.
\item
 $L^2_\mu$  is a Hilbert space equipped with the scalar product
 \b*\langle \delta, l\rangle_\mu := \int_{\R^*}\delta(u)l(u)\mu(du), \quad \text{ for all } \delta, l \in L^2_\mu\times L^2_\mu,\e*
 and the norm  $\|l\|^2_\mu:=\int_{\R^*}|l(u)|^2\mu(du)<\infty$.
\item    $\H^2$ is the set of real-valued predictable processes $\phi$ such that\\
 $\| \phi\|^2_{\H^2} := \E \left[\Int_0 ^T |\phi_t|^2 dt\right] < \infty.$
\item    $\H^2_\mu$ is the set of processes $l$ which are predictable, that is, measurable, \\
$l:(\Omega\times [0, T]\times U, \Pc\otimes \mathcal{U})\rightarrow (\R, \mathcal{B}(\R)); $ \quad $(\omega, t, u) \mapsto l_t(\omega, u)$ with \\
 $\| l\|^2_{\H^2_\mu} := \E \left[\Int_0 ^T \|l_t\|^2_\mu dt\right] < \infty.$
\item   ${\cal S}^2$ is the set of real-valued  optional  processes $\phi$ such that\\
$|||\phi|||^2_{{\cal S}^2}  := \E(\esssup_{\tau\in\stopo} |\phi_\tau |^2) <  \infty.$

\end{itemize}
\begin{Definition}[Driver, Lipschitz driver]\label{defd}
A function $f$ is said to be a {\em driver} if
\begin{itemize}
\item
$f: \Omega \times [0,T]\times \R^2 \times L^2_\mu \rightarrow \R $\\
$(\omega, t,y,z,k) \mapsto  f(\omega, t, y, z,k) $
  is $ {\Pc} \otimes {\cal B}(\R^2)\otimes {B}(L^2_\mu)$-measurable,
\item $E \left[(\int_0 ^T |f(t,0,0,0)| ^2 dt)\right]<\infty$.
\end{itemize}
A driver $f$ is called a {\em Lipschitz driver} if moreover there exists a constant $ K \geq 0$ such that $d\P \otimes dt$-a.e.\,,
for each $(y_1, z_1, k_1)$, $(y_2, z_2, k_2)$,
$$|f(\omega, t, y_1, z_1, k_1) - f(\omega, t, y_2, z_2, k_2)| \leq
K (|y_1 - y_2| + |z_1 - z_2|+ \|k_1 - k_2\|_\mu).$$
\end{Definition}
\begin{Assumption}
Assume that  $dP \otimes dt$-a.s\, for each $(y,z, k_1,k_2)$ $\in$ $ \mathbb{R}^2 \times (L^2_{\mu})^2$,
\be \label{assumpfjump}f( t,y,z, k_1)- f(t,y,z, k_2) \geq \langle \gamma_t^{y,z, k_1,k_2}  \,,\,k_1 - k_2 \rangle_\mu,\ee
with
\begin{equation*}
\gamma:  \Omega\times [0,T]  \times \mathbb{R}^2 \times  (L^2_{\mu})^2  \rightarrow  L^2_{\mu}\,; \, (\omega, t, y,z, k_1,k_2)\mapsto
\gamma_t^{y,z, k_1,k_2}(\omega,\cdot)
\end{equation*}
 ${\cal P } \otimes {\cal B}(\R^2) \otimes  {\cal B}( (L^2_{\mu})^2 )$-measurable, bounded, and satisfying $ dP\otimes dt \otimes \nu(du)$-a.s.\,, for each $(y,z, k_1,k_2)$ $\in$ $\mathbb{R}^2 \times (L^2_{\mu})^2$,
 \begin{equation}\label{condi}
\gamma_t^{y,z, k_1,k_2} (u)\geq -1\,\,\,  \;\; \text{ and }
\,\,  \;\;|\gamma_t^{y,z, k_1,k_2}(u)|  \leq \psi(u),
\end{equation}
where $\psi$ $\in$ $L^2_{\mu}$.\\
Hereafter, we denote for simplicity  $\gamma_t $ and we means $ \gamma_t^{y,z, k_1,k_2}$.
\end{Assumption}

Let $T>0$ be a fixed terminal time. Let $f$ be  a  driver.
Let $\xi= (\xi_t)_{t\in[0,T]}$ be in ${\cal S}^2$. 
We suppose moreover that the process $\xi$ is right USC. A process $\xi$ satisfying the previous properties will be called a  \emph{barrier}, or an  \emph{obstacle}.
\begin{Definition} \label{def_solution_RBSDE}
A process $(Y,Z,k, A,C)$ is said to be a solution to the reflected BSDE with parameters $(f,\xi)$, where $f$ is a driver and $\xi$ is a terminal value, if
\begin{align}\label{RBSDE}
&(Y,Z,k(\cdot),A,C)\in {\cal S}^2 \times \H^2\times \H^2_\nu \times {\cal S}^2\times {\cal S}^2 \nonumber \\
& Y_\tau=\xi_T+\int_{\tau}^T f(t,Y_t,  Z_t)dt-\int_{\tau}^T  Z_t dW_t -\int_{\tau}^T \int_U k_s(u)\tilde{N}(ds,du) +A_T-A_\tau+C_{T-} -C_{\tau-} \nn\\
& \text{ a.s. for all }\tau\in\stopo, \\
&  Y_t \geq \xi_t   \text{ for all } t\in[0,T]  \text{ a.s.,} \label{RBSDE_inequality_barrier}\\
& A \text{ is a nondecreasing right-continuous predictable process
with } A_0= 0, E(A_T)<\infty \text{ and such that } \nonumber\\
& \int_0^T {\bf 1}_{\{Y_t > \xi_t\}} dA^c_t = 0 \text{ a.s. and } \; (Y_{\tau-}-\bar{\xi}_\tau)(A^d_{\tau}-A^d_{\tau-})=0 \text{ a.s. for all (predictable) }\tau\in\stopo, \label{RBSDE_A}\\
& \text{where } \bar{\xi}_t:=  \limsup_{s \uparrow t, s< t} \xi_s. \nn\\
& C \text{ is a nondecreasing right-continuous adapted purely discontinuous process with } C_{0-}= 0, E(C_T)<\infty  \nonumber\\
& \text{ and such that }
(Y_{\tau}-\xi_{\tau})(C_{\tau}-C_{\tau-})=0 \text{ a.s. for all }\tau\in\stopo. \label{RBSDE_C}
\end{align}
\end{Definition}
We will use the following notation: Let $\beta >0$. For $\phi \in \H^2$, $\|\phi\|^2_\beta:=\E[\int_0^Te^{\beta s}\phi^2_sds]$. For $\phi \in \S^2$, $\vvertiii{\phi}^2_\beta:=\E[\esssup_{\tau\in \Tc_{0, T}}e^{\beta\tau}\phi^2_\tau]$. For $\phi \in \H^2_\mu$, $\|\phi\|^2_{\beta,\mu}:=\E[\int_0^T\int_Ue^{\beta s}|\phi_s(u)|^2\mu(du)ds]$.\\
In order to prove existence and uniqueness of the above reflected BSDE, we first investigate the case where the driver $f$ does not depend on $y,\,, z$ and $k$. In what follows, we adapt the proof of Grigorova et \textit{al}. \cite{GIOOQ15} to our setting. We first proof the a-priori estimate
\begin{Lemma}[A priori estimates]\label{Lemma_estimate}
Let $(Y^1, Z^1, k^1(\cdot), A^1, C^1)\in {\cal S}^2 \times \H^2 \times\H^2_\mu\times {\cal S}^2\times {\cal S}^2$ (resp. $(Y^2, Z^2,k^2(\cdot), A^2, C^2)\in {\cal S}^2 \times \H^2 \times\H^2_\mu\times {\cal S}^2\times {\cal S}^2$) be  a  solution to the RBSDE  associated with  obstacle $\xi$ and driver $f^1(\omega, t)$ (resp. $f^2(\omega, t)$).  There exists a positif constant $c>0$ such that for all $\varepsilon>0$,  for all $\beta\geq\frac 1 {\varepsilon^2}$ we have
\begin{equation}\label{eq_initial_Lemma_estimate}
\begin{aligned}
&\|Z^1-Z^2\|^2_\beta+\|k^1-k^2\|^2_{\beta,\mu}\,\leq \, \varepsilon^2  \|f^1-f^2\|^2_\beta\\
&\vvertiii{Y^1-Y^2}^2_\beta\,\leq \,2\varepsilon^2(1+4c^2)  \|f^1-f^2\|^2_\beta.
\end{aligned}
\end{equation}
\end{Lemma}
{\bf Proof.}
Let $\beta$ and  $\varepsilon$ be two nonnegative real numbers such that $\beta\geq\frac 1 {\varepsilon^2}$,
and   set \\
$\tilde Y:=Y^1-Y^2$, \\
$\tilde Z:=Z^1-Z^2$, \\
$\tilde k:=k^1-k^2$, \\
$\tilde A:=A^1-A^2$, \\
$\tilde C:=C^1-C^2$, \\
$\tilde f(\omega, t):=f^1(\omega, t)-f^2(\omega, t)$.\\
$\tilde Y_T=\xi_T-\xi_T=0$, \\
So that, we have
$$\tilde Y_\tau=\int_{\tau}^T \tilde f(t)dt-\int_{\tau}^T  \tilde Z_t dW_t-\int_{\tau}^T \int_U \tilde k_t(u)\tilde  N(du,dt) +\tilde A_T-\tilde A_\tau+\tilde C_{T-} -\tilde C_{\tau-} \text{ a.s. for all }\tau\in\stopo.$$
According to Theorem A.3 and Corollary A.2 in \cite{GIOOQ15}, and by noticing the decomposition $\tilde Y= \tilde Y_0+ M+A+B$, where
$M_t:=\int_0^t \tilde Z_s dW_s+\int_0^t \int_U \tilde k_s(u)\tilde  N(du,ds)$, $A_t:= -\int_0^t  \tilde f(s) ds- \tilde A_t$ and
$B_t:=-\tilde C_{t-}$, it is clear that  $\tilde Y$ is an optional (strong)  semimartingale.\\
Applying now \cite[Corollary A.2 ]{GIOOQ15} to $\tilde Y$ we obtain almost surely, for all $t\in[0,T]$,
\begin{equation}\label{eq0_lemma_estimate}
\begin{aligned}
 e^{\beta t}\tilde Y_t^2+\int_{]t,T]}  e^{\beta s} \tilde Z_s^2 ds+\int_{]t,T]} \int_U e^{\beta s} \tilde k_s^2(u) \mu(du)ds
 &\leq -\int_{]t,T]}\beta e^{\beta s} (\tilde Y_{s})^2 ds+2\int_{]t,T]}  e^{\beta s}\tilde Y_{s}\tilde f(s) ds\\
&\quad+2\int_{]t,T]}  e^{\beta s}\tilde Y_{s-}d\tilde A_s +2\int_{[t,T[}  e^{\beta s}\tilde Y_{s}d(\tilde C)_{s+}\\
&\quad-2\int_{]t,T]}  e^{\beta s}\tilde Y_{s-}\tilde Z_s d W_s\\
&\quad-2\int_{]t,T]} \int_U e^{\beta s}\tilde Y_{s-}\tilde k_s \tilde N(du,ds).
\end{aligned}
\end{equation}
From one side, we know that for any $a,b \in\R$ we have $2ab\leq (\frac a \varepsilon)^2+\varepsilon^2 b^2$, is follows a.e. for all $t\in[0,T]$,
\begin{equation*}
\begin{aligned}
-\int_{]t,T]}\beta e^{\beta s} (\tilde Y_{s})^2 ds+2\int_{]t,T]}   e^{\beta s}\tilde Y_{s}\tilde f(s) ds&
& \le(\frac 1 {\varepsilon^2}-\beta)\int_{]t,T]} e^{\beta s} (\tilde Y_{s})^2 ds+\varepsilon^2 \int_{]t,T]}   e^{\beta s}\tilde f^2(s) ds.
\end{aligned}
\end{equation*}
From other side, since $\beta\geq\frac 1 {\varepsilon^2}$, we have $(\frac 1 {\varepsilon^2}-\beta)\int_{]t,T]} e^{\beta s} (\tilde Y_{s})^2 ds \leq 0$, for all $t\in[0,T]$ a.s. \\
Now, following the same computation as in \cite{GIOOQ15} leads to a.s. for all $t\in[0,T]$
\b*
\int_{[t,T[}  e^{\beta s}\tilde Y_{s}d(\tilde C)_{s+}\leq 0 &\text{ and } &\int_{]t,T]}  e^{\beta s}\tilde Y_{s-}d\tilde A_s\le0 
\e*
The above observations, combined with equation \eqref{eq0_lemma_estimate}, lead to a.s. for all $t\in[0,T] $:
\begin{equation}\label{eq7_lemma_estimate}
\begin{aligned}
 e^{\beta t}\tilde Y_t^2+\int_{]t,T]}  e^{\beta s} \tilde Z_s^2 ds \leq \varepsilon^2 \int_{]t,T]}   e^{\beta s}\tilde f^2(s) ds
&-2\int_{]t,T]}  e^{\beta s}\tilde Y_{s-}\tilde Z_s d W_s-2\int_{]t,T]} \int_U e^{\beta s}\tilde Y_{s-}\tilde k_s \tilde N(du,ds),
\end{aligned}
\end{equation}
Next, we show that the two last terms on the r.h.s. of the previous inequality \eqref{eq7_lemma_estimate} have zero expectation. Using the left-continuity of a.e. trajectory of the process $(\tilde {Y}_{s-})$, we have
\begin{equation}\label{eq_a_lemma_estimate}
(\tilde{Y}_{s-})^2(\omega)\leq \sup_{t\in\Q} (\tilde Y_{t-})^2(\omega), \text{ for all } s\in(0,T],  \text{ for a.e. } \omega\in\Omega.
\end{equation}
On the other hand, for all $t\in (0,T]$, $(\tilde Y_{t-})^2\leq \esssup_{\tau\in\stopo} (\tilde Y_{\tau})^2$ a.s.; hence,
\begin{equation}\label{eq_b_lemma_estimate}
\sup_{t\in\Q}(\tilde Y_{t-})^2\leq \esssup_{\tau\in\stopo} (\tilde Y_{\tau})^2 \text{ a.s.}
\end{equation}
From equations \eqref{eq_a_lemma_estimate} and \eqref{eq_b_lemma_estimate}  together with Cauchy-Schwarz inequality, gives
\begin{equation*}
\begin{aligned}
E\left[ \sqrt{\int_{0}^T  e^{2\beta s}\tilde Y^2_{s-}\tilde Z^2_s ds}\right] &\leq
  E\left[\sqrt{\esssup_{\tau\in\stopo} \tilde (Y_{\tau})^2} \sqrt{\int_{0}^T  e^{2\beta s}\tilde Z^2_s ds}\right]
 \leq  \vvertiii{\tilde Y}_{{\cal S}^2} \|\tilde Z\|_{2\beta},
\end{aligned}
\end{equation*}
it follows that   $\vvertiii{\tilde Y}_{{\cal S}^2}<\infty$ since  $\vvertiii{\tilde Y}_{{\cal S}^2}\leq \vvertiii{Y^1}_{{\cal S}^2}+\vvertiii{Y^2}_{{\cal S}^2}$ and to $Y^1$ and $Y^2$ being in $\mathcal{S}^2$.
We also have that $ \|\tilde Z\|_{2\beta}<\infty$, due to the fact that $Z^1,Z^2\in\H^2$ and to the equivalence of the norms $\|\cdot\|_{2\beta}$ and $\|\cdot\|_{\H^2}$ (see the proof of Proposition 2.1 in \cite{GIOOQ15}).  We conclude that $E\left[ \sqrt{\int_{0}^T  e^{2\beta s}\tilde Y^2_{s-}\tilde Z^2_s ds}\right]<\infty$; whence, by standard arguments, we get
$E\left[ \int_{0}^T  e^{\beta s}\tilde Y_{s-}\tilde Z_s d W_s\right] =0.$

Similarly,
\begin{equation*}
\begin{aligned}
E\left[ \sqrt{\int_{0}^T \int_U e^{2\beta s}\tilde Y^2_{s-}\tilde k^2_s(u) \mu(du)ds}\right] &\leq
  E\left[\sqrt{\esssup_{\tau\in\stopo} \tilde (Y_{\tau})^2} \sqrt{\int_{0}^T  e^{2\beta s}\|\tilde k_s\|^2_\mu ds}\right]
 \leq  \vvertiii{\tilde Y}_{{\cal S}^2} \cdot \|\tilde k\|_{2\beta,\mu},
\end{aligned}
\end{equation*}
Since $ \|\tilde k\|_{2\beta,\mu}<\infty$, due to the fact that $k^1,k^2\in\H^2_\mu$ and to the equivalence of the norms $\|\cdot\|_{2\beta,\mu}$ and $\|\cdot\|_{\H^2_\mu}$.  We conclude that $E\left[ \sqrt{\int_{0}^T \int_U e^{2\beta s}\tilde Y^2_{s-}\tilde k^2_s(u) \mu(du)ds}\right]<\infty$; so that we conclude whence, by standard arguments, we get
$E\left[ \int_{]t,T]} \int_U e^{\beta s}\tilde Y_{s-}\tilde k_s (u)\tilde N(du,ds)\right] =0.$\\
Taking expectations on both sides of the inequality \eqref{eq7_lemma_estimate} with $t=0$, we get
$$\tilde Y_0^2+ \|\tilde Z\|^2_\beta+\|\tilde k\|^2_{\beta, \mu}\leq \varepsilon^2\|\tilde f\|^2_\beta.$$
Hence,
\begin{equation}\label{eq8_lemma_estimate}
\|\tilde Z\|^2_\beta+\|\tilde k\|^2_{\beta,\mu}\leq \varepsilon^2\|\tilde f\|^2_\beta.
\end{equation}
Again from the same inequality we have, for all $\tau\in\stopo$,
\begin{equation*}
\begin{aligned}
 e^{\beta \tau}\tilde Y_\tau^2\leq& \varepsilon^2 \int_{]0,T]}   e^{\beta s}\tilde f^2(s) ds
-2\int_{]0,T]}  e^{\beta s}\tilde Y_{s-}\tilde Z_s d W_s+2\int_{]0,\tau]}  e^{\beta s}\tilde Y_{s-}\tilde Z_s d W_s \\
&-2\int_{]0,T]} \int_U e^{\beta s}\tilde Y_{s-}\tilde k_s (u)\tilde N(du,ds)+2\int_{]0,\tau]} \int_U e^{\beta s}\tilde Y_{s-}\tilde k_s (u)\tilde N(du,ds) \text{ a.s.}
\end{aligned}
\end{equation*}
By taking first the essential supremum over $\tau\in\stopo$ and then the expectation on both sides of the above inequality, we obtain
\begin{equation}\label{eq9_lemma_estimate}
\begin{aligned}
E[\esssup_{\tau \in \stopo} e^{\beta \tau}\tilde Y_\tau^2]\leq &\varepsilon^2 \|\tilde f\|^2_\beta
+2E[\esssup_{\tau \in \stopo}|\int_0^\tau  e^{\beta s}\tilde Y_{s-}\tilde Z_s d W_s|]\\
&+2E[\esssup_{\tau \in \stopo}|\int_{]0,\tau]} \int_U e^{\beta s}\tilde Y_{s-}\tilde k_s (u)\tilde N(du,ds) |].
\end{aligned}
\end{equation}
By using the continuity of a.e. trajectory of the process $(\int_0^t  e^{\beta s}\tilde Y_{s-}\tilde Z_s d W_s)_{t\in[0,T]}$ and \\$(\int_{]0,\tau]} \int_U e^{\beta s}\tilde Y_{s-}\tilde k_s (u)\tilde N(du,ds) )_{t\in[0,T]}$ and Burkholder-Davis-Gundy inequalities (applied with $p=1$), we get
\begin{equation}\label{eq10_lemma_estimate}
2E[\esssup_{\tau \in \stopo}|\int_0^\tau  e^{\beta s}\tilde Y_{s-}\tilde Z_s d W_s|]=
2E[\sup_{t \in [0,T]}|\int_0^t  e^{\beta s}\tilde Y_{s-}\tilde Z_s d W_s|]\leq
2cE\left[ \sqrt{\int_{0}^T  e^{2\beta s}\tilde Y^2_{s-}\tilde Z^2_s ds}\right],
\end{equation}
and
\begin{equation}\label{eq11_lemma_estimate}
\begin{aligned}
2E[\esssup_{\tau \in \stopo}|\int_0^\tau  \int_U e^{\beta s}\tilde Y_{s-}\tilde k_s (u)\tilde N(du,ds)|]&=
2E[\sup_{t \in [0,T]}|\int_0^t  \int_U e^{\beta s}\tilde Y_{s-}\tilde k_s (u)\tilde N(du,ds)|]\\
&\leq
2cE\left[ \sqrt{\int_{0}^T\int_U  e^{2\beta s}\tilde Y^2_{s-}\tilde k^2_s (u)\mu(du)ds}\right],
\end{aligned}\end{equation}
where $c$ is a positive "universal" constant (which does not depend on the other parameters).
Following the same above reasoning we get:
$$\sqrt{\int_{0}^T  e^{2\beta s}\tilde Y^2_{s-}\tilde Z^2_s ds}\leq
\sqrt{\esssup_{\tau\in\stopo}  e^{\beta\tau}(\tilde Y_{\tau})^2\int_{0}^T  e^{\beta s}\tilde Z^2_s ds}\text{ a.s.} $$
and
$$\sqrt{\int_{0}^T\int_U  e^{2\beta s}\tilde Y^2_{s-}\tilde k^2_s \mu(du)ds}\leq
\sqrt{\esssup_{\tau\in\stopo}  e^{\beta\tau}(\tilde Y_{\tau})^2\int_{0}^T \int_U e^{\beta s}\tilde k^2_s(u) \mu(du)ds}\text{ a.s.}$$
Those inequalities with the the fact that  $ab\leq \frac 1 2 a^2+\frac 1 2 b^2$   leads to
$$2E[\esssup_{\tau \in \stopo}\int_0^\tau  e^{\beta s}\tilde Y_{s-}\tilde Z_s d W_s]\leq
\frac 1 2 E[\esssup_{\tau\in\stopo}  e^{\beta\tau}(\tilde Y_{\tau})^2]+ 2c^2E[\int_{0}^T  e^{\beta s}\tilde Z^2_s ds].$$
and
$$2E[\esssup_{\tau \in \stopo}|\int_0^\tau  \int_U e^{\beta s}\tilde Y_{s-}\tilde k_s (u)\tilde N(du,ds)|]\leq
\frac 1 2 E[\esssup_{\tau\in\stopo}  e^{\beta\tau}(\tilde Y_{\tau})^2]+ 2c^2E[\int_{0}^T  \int_U e^{\beta s}\tilde k^2_s(u)\mu(du) ds].$$
From this, together with \eqref{eq9_lemma_estimate}, we get
$$\frac 1 2 \vvertiii{\tilde Y}^2_\beta \leq \varepsilon^2 \|\tilde f\|^2_\beta+2c^2  \|\tilde Z\|^2_\beta+2c^2  \|\tilde k\|^2_{\mu,\beta}.$$
This inequality, combined with the estimate \eqref{eq8_lemma_estimate}  on $\|\tilde Z\|^2_\beta$ and $\|\tilde k\|^2_{\beta,\mu}$, gives
$$\vvertiii{\tilde Y}^2_\beta \leq 2\varepsilon^2(1+4c^2)  \|\tilde f\|^2_\beta.$$ \ep

We omit the proof of the following lemma since the proof is similar to \cite[Lemma 3.3]{GIOOQ15}
\begin{Lemma}\label{f}
Suppose that $f$ does not depend on $y,z$, that is $f(\omega,t, y,z) = f(\omega,t)$, where $f$
is a process in $\H^2$. Let $(\xi_t)$ be an obstacle. 
Then, the RBSDE from Definition \ref{def_solution_RBSDE} admits a unique solution $(Y,Z,k(\cdot),A,C)\in {\cal S}^2 \times \H^2\times \H^2_\mu \times {\cal S}^2\times {\cal S}^2$, and  for each $S \in \stopo$,  we have
\begin{eqnarray}\label{deux-2}
 Y_S= \esssup_{\tau \in\stops } E[ \xi_{\tau} + \int_S^\tau f(t)dt \mid \Fc_S] \quad   \rm{a.s.}
\end{eqnarray}
\end{Lemma}
We are going now to prove existence and uniqueness of RBSDE in Definition \ref{def_solution_RBSDE}, for general
driver using the fixed-point theorem
\begin{Theorem}\label{exiuni}
Let  $\xi$  be a left-limited and r.u.s.c. process in $\mathcal{S}^2$ and let $f$ be a  Lipschitz driver. \\
The RBSDE with parameters $(f,\xi)$ from Definition \ref{def_solution_RBSDE} admits a unique solution $(Y,Z,k(\cdot),A,C)\in \mathcal{S}^2  \times \H^2\times \H^2_\mu \times \mathcal{S}^2\times \mathcal{S}^2.$\\
Moreover, if $(\xi_t)$ is  assumed  l.u.s.c. along stopping times, then $(A_t)$ is continuous.
\end{Theorem}
{\bf Proof.}  \\
Let prove the uniqueness of the above RBSDE which will be obtained by via a fixed point of contraction of the function $\phi$ defined as follows:
Let:  $\mathcal{B}_\beta^2:=\mathcal{S}^2 \times \H^2\times \H^2_\mu$ endowed with the norm:
$$\| (Y, Z,k)\|_{\mathcal{B}_\beta^2}^2:=\vvertiii {Y}_{\beta}^2  +   \| Z\|_{\beta}^2+\| k\|_{\beta,\mu}^2, $$ for all $(Y,Z,k(\cdot))\in \mathcal{S}^2 \times \H^2\times \H^2_\mu.$\\
Let $\Phi$ be the map from $\mathcal{B}_\beta^2$ into itself which with $(Y, Z,k)$ associates $(Y, Z,k) = \Phi (y, z,\mathfrak{k})$, where $(y, z,\mathfrak{k})$ is the solution of RBSDE associated with with driver $f(s):=
f(s, y_s, z_s,\mathfrak{k}_s)$ and with obstacle $\xi$. Let $(A,C)$ be the associated Mertens process, constructed as in Lemma \ref{f}.  
The mapping $\Phi$ is well-defined by Lemma \ref{f}.  

%

%
Let $(y,z,k)$ and $(y',z',\mathfrak{k}')$ be two elements of $\mathcal{B}_\beta^2$. We set $(Y, Z,k)= \Phi (y, z,\mathfrak{k})$ and $(Y', Z',k') = \Phi (y', z'\mathfrak{k}').$  We also set $\tilde Y:=Y-Y'$, $\tilde Z:=Z-Z'$, $\tilde k:=k-k'$,  $\tilde y:=y-y'$ and $\tilde z:=z-z'$, $ \tilde{\mathfrak{k}}:=\mathfrak{k}-\mathfrak{k}'$.

According to Lemma \ref{Lemma_estimate}, we have
$$\vvertiii{\tilde Y}^2_\beta +\|\tilde Z\|^2_\beta+\|\tilde k\|^2_{\mu,\beta}\leq  \varepsilon^2(3+8c^2)  \|f(y,z,\mathfrak{k})-f(y',z',\mathfrak{k}')\|^2_\beta,  $$
for all $\varepsilon>0$,  for all $\beta\geq\frac 1 {\varepsilon^2}$. By using the Lipschitz property of $f$ and the convex property of square function we obtain
$$\|f(y,z,\mathfrak{k})-f(y',z',\mathfrak{k}')\|_{\beta}^2\leq C_{K}(\|\tilde y\|_{\beta}^2+\|\tilde z\|_{\beta}^2+\|\tilde{\mathfrak{k}}\|_{\beta,\mu}^2),$$
where $C_{K}$ is a positive constant depending on $K$. Which implies  for all $\varepsilon>0$,  for all $\beta\geq\frac 1 {\varepsilon^2}$, that
$$   \vvertiii{\tilde Y}^2_\beta +\|\tilde Z\|^2_\beta+\|\tilde k\|^2_{\beta,\mu}\leq  \varepsilon^2C_K(3+8c^2)(\|\tilde y\|_{\beta}^2+\|\tilde z\|_{\beta}^2+\|\tilde{\mathfrak{k}}\|_{\beta,\mu}^2).  $$
The previous inequality, combined with \cite[Remark 3.7]{GIOOQ15}, gives
$$   \vvertiii{\tilde Y}^2_\beta +\|\tilde Z\|^2_\beta+\|\tilde k\|^2_{\beta,\mu}\leq  \varepsilon^2C_K(3+8c^2)(T+1)(\vvertiii{\tilde y}_{\beta}^2+\|\tilde z\|_{\beta}^2+\|\tilde{\mathfrak{k}}\|_{\beta,\mu}^2).  $$
Now, Let  $\varepsilon>0$ be such that $\varepsilon^2C_K(3+8c^2)(T+1)<1$ and $\beta>0$ such that $\beta\geq\frac 1 {\varepsilon^2}$  the mapping $\Phi$ is a contraction on $\mathcal{B}_\beta^2$. Henceforth, there exist a triplet $(Y, Z, k) $  such that $\phi(Y, Z, k)=(Y, Z, k)$ which is the unique solution of the reflected RBSDE. 
\ep
\subsubsection{Dynamic risk measure induced by BSDEs with jumps}
In this subsection, we extend Proposition A.5 in \cite{GIOOQ15} to our setting, and we omit the proof of our Characterisation  Theorem, since its proof has the similar computation as in \cite[Theorem 4.2 ]{GIOOQ15} based mainly on Proposition \ref{compref}.\\
Let $T$ be a time horizon and $T'\in [0, T]$ a fixed instant before the terminal time $T$ . Let $f$ be a Lipschitz driver. Define the following functional: for each stopping time $\tau \in \Tc_{0, T}$ and $\xi \in \Sc^2$. Set
$$v(S)=-\esssup_{\tau \in \Tc_{S,T}} \mathcal{E}^f_{S,\t}(\xi_\tau), $$
where $S\in \Tc_{0, T}$,  $v$ is the dynamic risk measure, $\xi_{T'}$ is the gain of the position at time $T'$ and $-\mathcal{E}^f_{t,T'}(\xi_{T'})$ is the $f$-conditional expectation of $\xi_\t$ modelling the the risk at time $t$, where $t$ runs through interval $[0, T']$.
We are concerned to show that the minimal risk measure $v$ defined above coincide with $-Y$, where $Y$ is the solution of the reflected BSDE associated to $(f, \xi)$

\begin{Proposition}\label{compref} 
Let $f$ be a Lipschitz driver. Let $A$  be a nondecreasing right-continuous predictable process in ${\cal S}^2$ with $A_0=0$
 and let $C$ be a nondecreasing right-continuous adapted purely discontinuous process in ${\cal S}^2$ with $C_{0-}=0$.

 Let $(Y,Z,k) \in {\cal S}^2 \times \mathbb{H}^2\times \mathbb{H}^2_\mu$  satisfy
 \begin{equation}\label{dy}
  -d  Y_t  \displaystyle =  f(t,Y_{t}, Z_{t}, k_t)dt + dA_t + d C_{t-}-  Z_t dW_t-\int_Uk_t(u)\tilde{N}(du,dt),
  \end{equation}
 in the sense that, for each $\tau\in\stopo$, the equality
 \begin{eqnarray*}\label{dy}
 Y_{\tau}= Y_{T}+ \int_{\tau} ^T f(s,Y_{s}, Z_{s}, k_s)ds + A_T - A_{\tau} + C_{T-}- C_{\tau -}- \int_{\tau} ^T Z_{s}dW_{s}-\int_\tau^T\int_Uk_s(u)\tilde{N}(du,ds)
 \end{eqnarray*}
 holds almost-surely.
Then the process $(Y_t)$ is a strong ${\cal E}^f$-super-martingale (resp ${\cal E}^f$-submartingale).
\end{Proposition}
{\bf Proof.}
Let  $ \tau, \theta \in \T_0$ be such that $\tau \leq \theta$ a.s.\,
Let us show that $Y_{\tau} \geq {\cal E}^f_{\tau ,\theta}(Y_{\theta}) $ a.s.\\
We denote by $(X,\pi,l)$ the solution to the BSDE associated with driver $f$, terminal time $\theta$, and terminal condition $Y_{\theta}$; then ${\cal E}^f_{\tau ,\theta}(Y_{\theta})=  X_{\tau }$ a.s. (by definition of ${\cal E}^f$).
Set $\bar Y_t = Y_t - X_t$, $\bar Z_t = Z_t - \pi_t$, $\bar k_t(u) = k_t(u) - l_t(u)$. Then
  $$
  -d \bar Y_t  \displaystyle =  h_t dt + dA_t + d C_{t-}- \bar Z_t dW_t-\int_U\bar k_t(u)\tilde{N}(du,dt), \quad
   \bar Y_{\theta}  = 0,
   $$
where $h_t:=f (t, Y_{t-}, Z_t, k(\cdot)) - f(t, X_{t-}, \pi_t,l(\cdot))$.

By the same arguments as those of  the proof of the comparison theorem for BSDEs (c.f. \cite[Lemma 4.1]{QS13}), we have by assumption \eqref{assumpfjump}
\begin{equation}\label{eq2}
h_t \ge \varphi_t+\delta_t \bar {Y}_{t-} + \beta_t \bar Z_t+\langle\gamma_t,\bar k_t\rangle_\mu, \; 0<t \leq T,\quad
dP \otimes dt-{\rm a.e.,}\;
\end{equation}
where $\delta$ and $\beta$ are predictable bounded processes and $\gamma\ge-1$ $dt\otimes d\P\otimes \mu(du)$-a.s.

Let $\Gamma$ be  the unique solution to the following forward SDE with jumps
\begin{equation}\label{eq4b}
d \Gamma_{s}  = \displaystyle \Gamma_{s} \left[ \delta_s ds + \beta_s d W_s +\int_U \gamma_s(u)\tilde N(du,ds)\right], \;\;
\Gamma_{0}  = 1.
\end{equation}
Using the fact that $\bar Y$ is a strong optional semimartingale with decomposition  $\bar Y=M^1+A^1+B^1$, where $M^1_t=\int_0^t \bar Z_s dW_s+\int_0^t\int_U\tilde k_s(u)\tilde N(du,ds)$, $A^1_t:=-\int_0^t h_s ds-A_s$, and $B^1_t:=-C_{t-}$,
and applying the generalized change of variables formula from \cite[Theorem A.3]{GIOOQ15} with $n:=2$, $X^1:=\bar Y$, $X^2:= \Gamma$, and $F(x^1,x^2):=x^1x^2$. We obtain    
\begin{align*}
d(\bar Y_s \Gamma_s)
  = &  \Gamma_{s}(\bar Z_s+\bar Y_{s-}\beta_s)dW_s+\Gamma_s\int_U\big(\bar k_s(u)+\bar  Y_s\gamma_s(u)\big)\tilde N(du,ds)\\
  &+\Gamma\int_U\gamma_s(u)\bar k_s(u)N(du,ds)-\Gamma_s (dA_s + d C_{s}).
\end{align*}
Using inequality \eqref{eq2} together with the nonnegetivity of $\Gamma$ and doing the same computation as in \cite[Theorem 4.2 ]{QS13} we obtain
\begin{equation}\label{eq3}
  \bar Y_\tau\Gamma_\tau \ge E \left[  \int_\tau^\theta \Gamma_{s}\, (dA_s + d C_{s} )  \mid \Fc_{\tau} \right],\;\; \quad {\rm a.s.}
\end{equation}
Then, since $\Gamma_{s} \geq 0$, we have $\bar Y_\tau\Gamma_\tau \geq 0$ a.s. Since $\Gamma_\tau>0$ a.s., we have   $\bar Y_\tau\geq 0$, that is $Y_{\tau} \geq X_{\tau}={\cal E}_{\tau ,\theta}(Y_{\theta}) $ a.s.\,,The proof is thus complete.
\ep
\subsection{Refected BSDE with irregular barrier}
On a given complete probability space $(\Omega, \mathcal{F},\P)$, let $(B_t)_{0 \le t \le T }$ be a standard $d$-Brownian motion defined on a finite time interval $[0, T]$. Denote by $\mathbf{F}=\{\mathcal{F}_t\}_{0 \le t \le T } $ the augmentation of the natural filtration $\mathbf{F}^B=\{\mathcal{F}^B_t\}_{0 \le t \le T }$ with $\mathcal{F}^B_t:=\sigma\{B_s; 0 \le s \le t\}$ generated by $B$. We mean by $\mathcal{P}$ the $\sigma$-algebra of predictable set $\Omega\times [0, T]$.
For each $t\in [0, T]$ let introduce the following notations:
\begin{itemize}
    \item $\S^2$: the Banach space of all càdlàg process $Y$ such that
        \b*
        \|Y \|_{\S^2}:=\left(\E\left(\sup_{0\le t \le T}Y_t^2\right]   \right)
        \e*
    \item $D^2_\mathcal{F}:=\{v_t, 0 \le t \le T, \mbox{ is an }\mathcal{F}_t-\mbox{progressively measurable c�dl�g processes with } \\\E[\sup_{0 \le t \le T }|v_t|^2]<\infty\}$;
 \item $L^2_{\Fc}:=\{L:\Omega\times [0, T] \to \R, \mathbf{F}-\mbox{predictable bounded process}\}$;
  \item $\mathcal{H}^2:=\{v_t, 0 \le t \le T, \mbox{ is an }\mathcal{F}_t-\mbox{adapted process  such that }\E\int_0^T|v_t|^2dt<\infty\}$;
\end{itemize}
The generator $f: \Omega \times [0, T]\times \R^{1+d}\rightarrow\R^+ $ is a given $\mathcal{P}\bigotimes\mathcal{B}(\R^{1+d})$-measurable bounded function, and it satisfies:
\begin{itemize}
  \item The uniformly Lipschitz property w.r.t. $(y, z)$, i.e., there exists a constant $k\ge 0$ such that for any $y, y'\in \R$ $ z, z' \in \R^d$:
\begin{eqnarray}\label{f lipsh}
|f(\cdot,y,z)-f(\cdot,y',z')|&\leq& k\;\left(|y-y'|+|z-z' |\right)\s \P-\mbox{a.s..}\end{eqnarray}
  \item $  \E[\Int_0^T f(\omega, t, 0, 0)^2 dt]<\infty.$
\end{itemize}
Consider the following reflected BSDE with $[\mathfrak{B}]$--obstacle associated to $(f, L)$. A solution is a triplet  $(Y, Z, K):=(Y_t, Z_t, K_t)_{t\le T}$ of processes with value in $\R^{1+d}\times\R^+$, which satisfies :
\begin{equation}\label{RBSDET}
\left\{
\begin{array}{lll}
Y_t = -\Int_0^tf(s, Y_s, Z_s)ds+\Int_t^T dK_s- \int_t^T Z_sdB_s; \\
Y\geq L \s \mbox{ a.s. a.e.};\\
\Int_0^{T}(Y_{s^-} - L^*_{s^-})dK_s=0, \s \forall L^* \in D_\mathcal{F}\mbox{ such that } L \leq L^* \leq Y \s \mbox{a.s. a.e.}
\end{array}
\right.
\end{equation}
It is worth mentioning here, that a similar kind of RBSDE was first studied in Peng and Xu \cite{PX05} using penalization method and in Essaky et al. \cite{EHO12} using generalized Snell envelop. And recently, a more general setting was studied by Grigorova et al. \cite{GIOOQ15} using Strong supermartinglae and the theory of optional process. First, let us assume the map $f$ does not depend on $(y, z)$, i.e., $\P$-a.s., $f(\omega, t, y, z)=g(\omega , t)$, for any $t, y$ and $z$. In the following we establish the existence of the solution of the RBSDE associated to $(g, L)$.
\begin{Theorem} \label{RBSDEone}
We assume that the lower obstacle $L \in [\mathfrak{B}] $. Then the reflected BSDE \eqref{RBSDET} has a unique minimal solution $(Y,Z,K) \in D^2_\mathcal{F} \times \mathcal{H}^2\times D^2_\mathcal{F}$ associated to $(g, L)$. Moreover, $\Big(Y_t + \Int_0^tg(s)ds\Big)_{0 \le t \le T }$ is the Strong envelope of $\Big(L_t + \Int_0^tg(s)ds\Big)_{0 \le t \le T }$.
\end{Theorem}
{\bf Proof.}  \\
Let $\bar{Y}$ be the Strong envelope of $\bar{L}:=L+\int_0^\cdot g(s)ds $, then $dt\otimes d\P $-a.e.
$\bar{Y}\ge \bar{L}. $ By Doob-Meyer decomposition of the super-martingale $\bar{Y}$, there exist an increasing c�dl�g process $\bar{K}$ with $\bar{K}_0=0$ and a continuous uniformly integrable martingale $\bar{M}$ such that  $\bar{Y}_t=\bar{M}_t - \bar{K}_t$.  From Proposition \ref{PropSE(X)}, for any $\bar{L}^* \in \S^2$ such that $\bar{Y}\ge \bar{L}^* \ge \bar{L} $ a.s. a.e. we have:
\begin{equation}\label{SkorohodKbar}
\int_0^T (\bar{Y}_{s^-}-\bar{L}^*_{s^-} )d\bar{K}_s=0. \end{equation}
Denoting $Y:=\bar{Y}-\int_0^\cdot g(s)ds$, it follows that
$$Y \ge L^*:=\bar{L}^* -\int_0^\cdot g(s)ds \ge L, \s \mbox{ a.s. a.e.}$$
From other side, for any $\bar{L}^* \in \S^2$, we have $L^* \in \S^2$, this with \eqref{SkorohodKbar} leads to
$$\int_0^T (Y_{s^-}-L^*_{s^-} )d\bar{K}_s=0. $$
Since $(\bar{Y}_t)_{0 \le t \le T }$ is a square integrable super-martingale, thanks to Remark \ref{U in S2}. we have $\E [\bar{K}^2_T] <\infty$   (Dellacherie and Meyer  \cite{DM82}). Hence the following martingale
 \begin{eqnarray*}
 \bar{M}_t=\E\left[\bar{M}_T/\mathcal{F}_t\right]=\E\left[\bar{K}_T/\mathcal{F}_t\right],
 \end{eqnarray*}
 is also a square integrable, and by martingale representation theorem, there exists a unique predictable process $Z$ such that
 \begin{eqnarray*}
 \bar{M}_t=\bar{K}_T-\int_t^T Z_sdB_s,
 \end{eqnarray*}
where $\E\int_0^T |Z_t|^2dt<\infty$, i.e. $Z\in \mathcal{H}^2$. It follows that
\begin{eqnarray*}
Y_t &=&-\int_0^{t}g(s)ds - \int_t^{T} Z_sdB_s+\bar{K}_T- \bar{K}_t\\
    &\geq& L_t, \s dt\otimes d\P-\mbox{a.e.}
\end{eqnarray*}
Hence, the existence of solution  $(Y, Z, K)$ for RBSDE (\ref{RBSDET}).
\\
Now, we return to show uniqueness of the solution. We suppose the existence of  two solutions $(Y, Z, K)$ and $(Y', Z', K')$ of RBSDE (\ref{RBSDET}), s.t $Y' \geq L$ a.s. a.e. and $Y_T=Y'_T=0$ a.s.
Then
\begin{eqnarray*}
Y'_t+ \int_0^{t}g(s)ds &=& \mathbb{E}\left[ K'_T / \mathcal{F}_t\right] - K'_t \\
                    &\geq& L_t + \int_0^{t}g(s)ds, \s dt\otimes d\P -\mbox{a.s.}.
\end{eqnarray*}
According to Proposition \ref{PropSE(X)}, we have
\begin{eqnarray*}
Y'_t+ \int_0^{t}g(s)ds &\geq& SE\left(L_t + \int_0^{t}g(s)ds \right)\\
                        &\geq& Y_t + \int_0^{t}g(s)ds, \s\P-\mbox{a.s.},
\end{eqnarray*}
which leads to the minimality of the solution.
\ep\\
We are now ready to give the main result of this section
\begin{Theorem}
The reflected BSDE \eqref{RBSDET} associated to $(f, L)$ has a unique minimal solution $(Y, Z, K)$.
\end{Theorem}
{\bf Proof. }
We show existence, which will be obtained via a fixed point of the contraction of the function $l$ defined  as follows:\\
Let $\Dc:= D^2_\mathcal{F}\times \Hc^2$ endowed with the norm,
\begin{eqnarray}
\|(Y, Z)\|_\beta = \left\{\E\left[\int_0^Te^{\beta s}|Y_s|^2+e^{\beta s}|Z_s|^2ds  \right]\right\}^\frac12; \s \beta >0.
\end{eqnarray}
Let $l$ be the map from $\Dc$ to itself   which with $(Y, Z)$ associates $l(Y, Z)=(\bar{Y}, \bar{Z})$ where $(\bar{Y}, \bar{Z})$  is the solution of the reflected BSDE associated with $(f(\cdot, Y, Z), L)$. Let $(Y', Z')$  be another couple  of $\Dc$ and $l(Y,' Z')=(\bar{Y}', \bar{Z}')$, then using It�'s formula we obtain, for any $t\le T$,
\begin{eqnarray*}
&&e^{\beta t}(\Yb_t-\Yb_t')^2+\beta\int_t^T e^{\beta s}(\Yb_s-\Yb_s')^2ds+\int_t^T e^{\beta s}(\Zb_s-\Zb_s')^2ds\\
&&=(M_T -M_t)+2\int_t^T e^{\beta s}(\Yb_s-\Yb_s')(d\bar{K}_s-d\bar{K}'_s)\\
&&\s-2\int_0^t e^{ \beta s}(\Yb_s-\Yb_s')(f(s, Y_s, Z_s)-f(s, Y_s', Z_s'))ds
\end{eqnarray*}
where $M$ is the martingale part. We set $L^*_t=\Yb_t\wedge\Yb_t'$, is follows that $L^* \in D^2_\mathcal{F}$ satisfies $L_t\le L_t^* \le \Yb $ and $L_t\le L_t^* \le \Yb' $. So that the general Skorohod condition is satisfied
$$\Int_0^{T }(\Yb_{s^-} - L^*_{s^-})d\bar{K}_s=\Int_0^{T }(\Yb_{s^-}' - L^*_{s^-})d\bar{K}_s'=0.$$
Consequently, $\int_t^T e^{\beta s}(\Yb_s-\Yb_s')(d\bar{K}_s-d\bar{K}'_s) \le 0$. Then
\begin{eqnarray*}
&&\beta \E[\int_0^T e^{\beta s}(\Yb_s-\Yb_s')^2ds]+\E[\int_0^T e^{\beta s}(\Zb_s-\Zb_s')^2ds]\\
&& \s\le2\E[\int_0^T e^{\beta s}|\Yb_s-\Yb_s'|\cdot|f(s, Y_s, Z_s)-f(s, Y_s', Z_s')|ds]\\
&&\s \le k\eps\E[\int_0^T e^{\beta s}(\Yb_s-\Yb_s')^2ds]+\frac{k}{\eps}\E[\int_0^T e^{\beta s }\{|Y_s-Y_s'|^2+|Z_s-Z_s'|^2\}ds].
\end{eqnarray*}
Hence
\begin{eqnarray*}
&&(\beta-k\eps)\E[\int_0^T e^{\beta s}(\Yb_s-\Yb_s')^2ds]+\E[\int_0^T e^{\beta s}(\Zb_s-\Zb_s')^2ds]\\
&&\s \le \frac{k}{\eps}\E[\int_0^T e^{\beta s }\{|Y_s-Y_s'|^2+|Z_s-Z_s'|^2\}ds].
\end{eqnarray*}
For $\beta$ bigger enough and $\eps$ such that $k\le \eps \le \frac{\beta-1}{k} $, then $l$ is a contraction on $\Dc$, so that ther exists a couple $(Y, Z)$ such that $l(Y, Z)=(Y, Z)$ which, with $K$, is the unique solution of the reflected BSDE associated to $(f, L)$.
\ep

\end{document}